\documentclass[12pt]{article}
\usepackage{fancyhdr, array,calc,graphicx,url,tabularx}
\usepackage{geometry}
\usepackage{latexsym}
\usepackage{amssymb}
\usepackage{amsmath}
\usepackage{makeidx}
\usepackage{enumerate}
\usepackage{verbatim}
\usepackage{tikz}
\usepackage[colorlinks=true,linkcolor=blue]{hyperref}

\usepackage{changepage}

\makeindex

{
   \end{minipage}
   \vspace*{\stretch{3}}
   \clearpage
}

\renewcommand{\gg}{\gamma}

\newcommand{\bR}{{\bf P}}

\newcommand{\rest}{\restriction}

\newcommand{\forces}{\Vdash}

\renewcommand{\models}{\vDash}
\newcommand{\powerset}{{\cal P}}

\newcommand{\insegeq}{\trianglelefteq}

\newcommand{\cp}{{\rm cp }}

\newcommand{\K}{{ k:V \to V }}

\newtheorem{theorem}{Theorem}[section]

\newtheorem{definition}[theorem]{Definition}

\newtheorem{lemma}[theorem]{Lemma}
\newtheorem{corollary}[theorem]{Corollary}

\newtheorem{question}[theorem]{Question}

\numberwithin{figure}{section}

\newenvironment{proof}{{\it{
Proof.}}}{\nopagebreak\mbox{}{\hfill$\square$}
\par\bigskip}
\newenvironment{Claim}{\noindent{\it
Claim} {}}{}
\newenvironment{Case}{\noindent{\it
Case}}{}

\newcommand{\rthm}[1]{Theorem~\ref{#1}}
\newcommand{\rlem}[1]{Lemma~\ref{#1}}

\newcommand{\rdef}[1]{Definition~\ref{#1}}

\newcommand{\rques}[1]{Question~\ref{#1}}

\def\inseg{\triangleleft}

\def\k{\kappa}
\def\a{\alpha}
\def\b{\beta}
\def\d{\delta}

\def\l{\lambda}

\def\P{{\mathcal{P} }}
\def\W{{\mathcal{W} }}
\def\Q{{\mathcal{ Q}}}

\def\K{{\mathcal{ K}}}

\def\R{{\mathcal R}}

\def\H{{\rm{HOD}}}
\def\M{{\mathcal{M}}}
\def\N{{\mathcal{N}}}

\def\T {{\mathcal{T}}}
\def\U{{\mathcal{U}}}
\def\S{{\mathcal{S}}}

\def\bR{{\mathbb{R}}}

\def\VT{{\vec{\mathcal{T}}}}

\def\cp #1{{ crit  #1 }}

\def\iff{\mathrel{\leftrightarrow}}

\def\and{\mathrel{\kern1pt\&\kern1pt}}

\def\<#1>{\langle\,#1\,\rangle}

 \input xy
 \xyoption{all}

\title{Translation procedures in descriptive inner model theory
}

\author{Grigor Sargsyan\thanks{The author's work is partially based upon work supported by the National Science Foundation under Grant No DMS-1352034 and DMS-1201348.}\\
        Department of Mathematics\\
        Rutgers University\\
        100 Frelinghuysen Rd.\\
        Piscataway, NJ 08854 USA\\
        http://math.rutgers.edu/$\sim$gs481\\
        grigor@math.rutgers.edu\\\\
        }   

\date{\today}

\begin{document}

\maketitle

\begin{abstract}
We develop a basic translation procedure that translates a hod mouse to an equivalent mouse. Unlike the translation procedure used by Steel and Zhu, our procedure works in a coarse setting without assuming $AD^+$. Nevertheless, the procedure resembles the one developed by Steel in \cite{DMATM}. We use the translation procedures to answer a question of Trevor Wilson. Namely, we show that if there is a stationary class of $\l$ such that $\l$ is a limit of Woodin cardinals and the derived model at $\l$ satisfies $AD^++\theta_0<\Theta$ then there is a transitive model $M$ such that $Ord\subseteq M$ and $M\models ``$there is a proper class of Woodin cardinals and a strong cardinal". Using a theorem of Woodin on derived models it is not hard to see that the reverse of the aforementioned theorem is also true, thus proving that the two theories are in fact equiconsistent. 
\end{abstract}

In descriptive inner model theory, the role of hod mice is to capture descriptive set theoretic strength of sets of reals inside a fine structural model. Such capturing leads to representation of $\H$ of models of $AD^++V=L(\powerset(\bR))$ as a fine structural model (see Theorem 4.24 of \cite{ATHM}). It also leads to partial proofs of the \textit{Mouse Set Conjecture} (see Theorem 6.19 of \cite{ATHM}). For more on the role of hod mice in descriptive inner model theory see \cite{BSL}. 

While the theory of hod mice has been advanced to the level of the theory \textit{Largest Suslin Axiom} (\textit{LSA})\footnote{LSA is the conjunction of the following axioms, $AD^+$, $\neg AD_{\bR}$, $V=L(\powerset(\bR))$ and ``the largest Suslin cardinal is a member of the Solovay sequence", see \cite{hodlsa} and \cite{complsa}.} the relationship between hod mice and ordinary mice is not as well understood. Steel, in \cite{DMATM}, initiated the study of the relationship between hod mice and ordinary mice. Yizheng Zhu, in unpublished work, building on ideas of Steel, has shown that many models of $AD^++V=L(\powerset(\bR))$ are derived models of mice. He devised a translation procedure that \textit{translates} a hod mouse, which is below the theory $AD_{\bR}+``\Theta$ is regular", into a mouse while making sure that the derived models of the original hod mouse and the resulting mouse are the same. What is being \textit{translated} in these constructions is the strategy predicate of the hod mouse in question. It is being translated into an ordinary extender sequence.

The motivating question that has started the current project is whether such translation can be done locally. In both Steel's work and Zhu's work, the translation is done using $AD^+$ theory. Their procedure doesn't work inside an arbitrary hod mouse but rather a specially chosen one. In this paper, we devise a translation procedure that can be done locally. It eliminates several complications of the two aforementioned constructions. We will then use our construction to answer a question of Wilson (see \rques{trev question}).\\\\

\textbf{Acknowledgments.} The author would like to thank Trevor Wilson for asking \rques{trev question}. He also wishes to thank the referee for a very long list of improvements. Unfortunately time constraints made it impossible for the author to address every invaluable suggestion made by the referee. 

\section{Translatable structures}

We start by introducing the hybrid models that we will translate into ordinary mice. We note once again that the general structure of our construction follows that of Steel's construction introduced in Section 15 of \cite{DMATM}. The objects introduced in this section are quite standard and their prototypes first appeared in Woodin's work on the analysis of $\H$ (see, for instance, Section 3 of \cite{StrengthPFA1} or Section 1 of \cite{TrangSarg}). 

Before we start, we remark that definition of the premouse in this paper is somewhat different than in other papers in print. We will use an indexing scheme in which extenders with critical point $\k$, which is a limit of Woodin cardinals and is a cutpoint, are indexed at the successor of the least cutpoint of their own ultrapower. More precisely, if $\M$ is a premouse and $E_\a \in \vec{E}^\M$ is an extender then if $\cp(E_\a)$ is a cutpoint in $\M$ and $\M\models ``\cp(E_\a)$ is a limit of Woodin cardinals" then in $Ult(\M, E)$, $\a$ is the successor of $\nu$ where $\nu\geq \k^+$ is the least cutpoint of $Ult(\M, E)$.

Suppose that $V\models ZFC+$``there is a proper class of Woodin cardinals". We say that the transitive set $X$ is \textit{self-well-ordered} if there is a well-ordering of $X$ in $\mathcal{J}_\omega(X)$. We then let 

\begin{definition}\label{stack}
\begin{center}
$\W(X)=\cup\{ \M: \M$ is an $Ord$-iterable sound $X$-mouse such that $\rho_\omega(\M)=X\}$.
\end{center}
We then, for $\a\in Ord$, define $\W_\a(X)$ by induction using the following scheme:
\begin{enumerate}
\item $\W_0(X)=\W(X)$.
\item $\W_{\xi+1}=\W(\W_\xi(X))$.
\item For limit ordinal $\l$, $\W_\l(X)=\cup_{\xi<\l}\W_\xi(X)$.
\end{enumerate}
\end{definition}

Suppose now that $\P$ is a premouse. We then say 

\begin{definition}[Suitable Premouse]\label{suitable}
$\P$ is suitable if there is an ordinal $\d$ such that
\begin{enumerate}
\item $\P\models ``\d$ is the unique Woodin cardinal",
\item for every strong $\P$-cutpoint\footnote{We say $\xi$ is $\P$-cutpoint if there is no $E\in \vec{E}^\P$ such that $\xi\not \in (\cp(E), \a(E))$ where $\a(E)$ is the index of $E$. We say $\eta$ is a strong cutpoint if for any $E$ as above, $\xi\not \in[\cp(E), \a(E))$.} $\eta$, $\P|(\eta^+)^\P=\W(\P|\eta)$. 
\item $\P=\W_\omega(\P|\d)$.
\end{enumerate}
\end{definition}

Suppose now that $\P$ is a suitable premouse. We let $\d^\P$ be the Woodin cardinal of $\P$. Given an iteration strategy $\Sigma$ for $\P$, we say $\Sigma$ is \textit{fullness preserving} if whenever $\Q$ is a $\Sigma$-iterate of $\P$ via $\Sigma$ then $\Q$ is suitable (see, for instance, Definition 6.2 of \cite{StrengthPFA1}). We say $\Sigma$ has \textit{branch condensation} if whenever $i: \P\rightarrow \Q$ is an iteration via $\Sigma$, $\VT$ is a stack according to $\Sigma$ such that $\VT\in dom(\Sigma)$, $b$ is a branch of $\VT$ such that $\pi^{\VT}_b$ exists and there is $\sigma:\M^{\VT}_b\rightarrow \Q$ with the property that $i=\sigma\circ \pi^{\VT}_b$, then $b=\Sigma(\VT)$ (see Definition 2.14 of \cite{ATHM}). Both fullness preservation and branch condensation are important inner model theoretic notions and have been fully explored in literature (see \cite{ATHM}). 

We say $(\P, \Sigma)$ is a hod pair if $\P$ is suitable and $\Sigma$ is an $(<Ord, <Ord)$-iteration strategy for $\P$ that is fullness preserving and has branch condensation. Notice that all the notions introduced thus far are relative to $M$, and to emphasize this we will sometimes say ``$M$-suitable", ``$M$-fullness preserving" and etc. Next we introduce $M$-super fullness preservation, which also comes from \cite{ATHM} (see Definition 3.33 of \cite{ATHM}). This notion is key in showing that the strategy can be interpreted on to generic extensions (see, for instance, Lemma 3.35 of \cite{ATHM}). We continue working in $M$. 

\begin{definition}[Super fullness preservation]\label{super fullness preservation}
Suppose $(\P, \Sigma)$ is a hod pair such that $\Sigma$ has branch condensation and is fullness preserving. Let $\d=\d^\P$ and for $n\in \omega$, let $\d_n=(\d^{+n})^\P$. We say $\Sigma$ is $\eta$-super fullness preserving if for every $n<\omega$, there is a term $\tau\in \P^{Coll(\omega, (\d^{+n})^\P)}$ such that 
\begin{enumerate}
\item $\P\models \forces_{Coll(\omega, \d_n)}\tau\subseteq \bR^2$,
\item whenever $g\subseteq Coll(\omega, \d_n)$ is $\P$-generic, 
\begin{center}
$\tau_g\subseteq \{(x, y): x\in \bR \wedge y$ codes a sound $x$-mouse projecting to $\omega\}$.
\end{center}
\item whenever $j:\P\rightarrow \Q$ is an iteration according to $\Sigma\rest V_\eta^M$, $X\in V_\eta^M$ is a self-well-ordered set and $g\subseteq  Coll(\omega, j(\d_n))$ is $\Q$-generic such that $X\in HC^{\Q[g]}$, letting $x\in \bR^{\Q[g]}$ be such that $x$ codes $X$ and $x$ is generic over $\mathcal{J}_\omega[X]$,
\begin{center}
$\W(X)=\cup\{ \M: \exists y ( (x, y)\in (j(\tau))_g \wedge y$ codes the $\S$-translate of $\M$ to an $x$-mouse$\}$.
\end{center}
\end{enumerate}
\end{definition}

Recall that using $\S$-constructions, it is possible to translate  an $X$-mouse into an $x$-mouse and vice versa. See Section 3.8 of \cite{ATHM} for more on $\S$-translations. We need one more notion before we can introduce translatable structures.

Suppose $\Sigma$ is super fullness preserving. We say $M$ is \textit{generically closed under} $\Sigma$ if in all generic extensions of $M$, $\Sigma$ has super fullness preserving extension. More precisely, 

\begin{definition} We say $M$ is generically closed under $\Sigma$ if for all $M$-cardinals $\eta$ and for all $\leq \eta$-generics $g$ there is a unique $\Lambda\in M[g]$, denoted by $\Sigma^{g, \eta}$, such that $\Lambda$ is an $((\eta^+)^M, (\eta^+)^M)$-strategy, $\Lambda\rest M$ is the $((\eta^+)^M, (\eta^+)^M)$-fragment of $\Sigma$ and $M[g]\models ``\Lambda$ is super fullness preserving".
\end{definition}

Suppose now that $M$ is generically closed under $\Sigma$. Given an $M$-generic $g$, we let $\Sigma^g=\cup_{\eta\in Ord} \Sigma^{g, \eta}$. 

Suppose next that $\P$ is suitable. Let $\d^\P_n=((\d^\P)^{+n})^\P$. Suppose $\tau\in \P^{Coll(\omega, \d_n^\P)}$ is a term relation for a set of reals.  We then let 
 \begin{center}
 $\gg_\tau^\P=\sup (Hull_1^\P(\{ \tau\})\cap \d^\P)$ and\\
 $H^\P_\tau=Hull_1^\P(\{\tau\}\cup \gg_\tau^\P)$.
 \end{center}

Suppose now that $\Sigma$ is a super fullness preserving strategy for $\P$ such that $M$ is generically closed under $\Sigma$.  We then say 

\begin{definition} $\Sigma$ guides $\tau$ correctly if for any $M$-generic $g$, $\Sigma^g$-iterates $\Q$ and $\R$ of $\P$ such that the iteration embeddings $i:\P\rightarrow \Q$ and $j:\P\rightarrow \R$ exist, and an integer $n<\omega$, whenever $x\in \mathbb{R}^{M[g]}$, $k\subseteq Coll(\omega, i(\d_n))$ and $l\subseteq Coll(\omega, j(\d_n))$ are such that $k$ is $\Q$-generic, $l$ is $\R$-generic and $x\in \Q[k]\cap \R[l]$ then 
\begin{center}
 $x\in i(\tau)_k\iff x\in j(\tau)_l$
 \end{center}
 \end{definition}
 We let $A_{\Sigma, \tau}^g$ be the set of reals in $M[g]$ determined by $\Sigma^g$ and $\tau$. Notice that for any $\Q$ as above $i(\tau)$ \textit{term captures} $A_{\Sigma, \tau}^g$ (see Definition 2.18 of \cite{ATHM}). 
 
Next we identify strategies that are guided by term relations (see Definition 4.14 of \cite{ATHM}).

\begin{definition}\label{guided by term relations} Suppose $(\P, \Sigma)$ is a hod pair such that $\Sigma$ is super fullness preserving and $M$ is generically closed under $\Sigma$. We say $\Sigma$ is guided by the sequence $(\tau_i: i<\omega)$ if for any $M$-generic $g$
\begin{enumerate}
\item for every $i$ there is $n$ such that $\tau_i\in \P^{Coll(\omega, \d_n^\P)}$ is a term relation for a set of reals,
\item for each $i<\omega$, $\Sigma^g$ guides $\tau_i$ correctly,
\item given any $\Sigma^g$-iterate $\Q$ of $\P$ such that the iteration embedding $\pi:\P\rightarrow \Q$ exists, $\Q=\cup_{i<\omega}H^\Q_{\pi(\tau_i)}$.
\item given any $\Sigma^g$ iterate $\Q$ of $\P$ such that the iteration embedding $\pi:\P\rightarrow \Q$ exists and given any normal tree $\T\in dom(\Sigma^g_\Q)$, letting $b=\Sigma^g_\Q(\T)$, if $\pi^\T_b$ exists then $b$ is the unique branch of $\T$ such that for every $i<\omega$, $\pi^\T_b(\pi(\tau_i))$ locally term captures $A^g_{\Sigma, \tau_i}$,
\item given any $\eta$ that is a limit of Woodin cardinals of $M[g]$ such that $\P\in V_\eta^{M[g]}$, if $h\subseteq Coll(\omega, <\eta)$-generic then for every $i<\omega$, $A^{g*h}_{\Sigma, \tau_i}$ is an $OD$ set of reals in the derived model of $M[g]$ as computed by $h$.
\end{enumerate}
We say $\Sigma$ is strongly guided by the sequence $(\tau_i: i<\omega)$ if whenever $\pi:\P\rightarrow \Q$ comes from an iteration according to $\Sigma^g$ and $\R$ is such that there are $\sigma:\P\rightarrow \R$ and $k:\R\rightarrow \Q$ such that $\pi=k\circ \sigma$ then $\R$ is suitable and for every $i$, $k^{-1}(\pi(\tau_i))$ term captures $A^g_{\Sigma, \tau_i}$.
\end{definition}

We can now introduce the structures that we will translate into ordinary premice. 

\begin{definition}[Translatable structure]\label{translatable structure}
We say $\K$ is a translatable structure if $\K$ is a hybrid strategy premouse (see \cite[Section 1.2]{ATHM}) satisfying the following conditions.
\begin{enumerate}
\item $\K\models ``$there are proper class of Woodin cardinals and there is no measurable limit of Woodin cardinals".
\item Let $\d$ be the least Woodin cardinal of $\K$. In $\K$, let $\P=\W_\omega(\K|\d)$. Then $\P$ is $\K$-suitable, and 
$\Sigma^\K$ is a $\K$-super fullness preserving iteration strategy for $\P$ such that $\K$ is generically closed under $\Sigma^\K$ 
and there is a collection of term relations $A\subseteq \P$ such that whenever $g\subseteq Coll(\omega, \P)$ is $\K$-generic and $(\tau_i: i<\omega)$ is a generic enumeration of $A$ in $\K[g]$, $(\Sigma^\K)^g$ is strongly guided by $(\tau_i : i<\omega)$. 
\item Suppose $\kappa<\lambda$ are two consecutive Woodin cardinals greater than $\d$. Let $\nu\in (\k, \l)$ be a cardinal of $\K$. Then for any $\eta$, $\K|(\nu^+)^\K$, as a $\Sigma^\K$-premouse, has an $(\eta, \eta)$-iteration strategy $\Lambda\in \K$ for stacks that are based on the window $(\k, \nu)$ such that $\Lambda$ is also a strategy for $\K$\footnote{Note that we easily get that $\K$ is internally $(Ord, Ord)$-iterable for stacks based on the window $(\k, \nu)$.}. 
\item $\K\models ``$there is no inner model with a superstrong cardinal".
\item Suppose $(\d_i: i<\omega)$ are the first $\omega$ many Woodin cardinals of $\K$ and $\d_\omega$ is their limit. Suppose $g\subseteq Coll(\omega, <\d_\omega)$ is $\K$-generic and $M$ is the derived model computed in $\K(\bR^*)$\footnote{Thus, $M=\mathcal{J}(Hom^*, \bR^*)$ where $\bR^*=\cup_{\xi<\d_\omega}\bR^{\K[g\cap Coll(\omega, \xi)}$ and in $\K(\bR^*)$, $Hom^*=\{ A\subseteq \bR^*:$ for some $\xi<\d_\omega$ there is a pair of $\d_\omega$-absoluteley complementing trees $(S, T)\in \K[g\cap Coll(\omega, \xi)]$ such that $(p[S])^{\K(\bR^*)}=A\}$.}. Let $\Phi=(\Sigma^\K)^g\rest HC^M$.Then $\Phi\in M$ and in $M$, $\K|(\d_1^{+\omega})^\K$, as a $\Phi$-premouse, is suitable and $\emptyset$-iterable (for the definitions of these notions see Section 1 of \cite{TrangSarg})\footnote{$\emptyset$-iterable refers to the following. Suppose $\P=\K|(\d_1^{+\omega})^\K$ and $\T$ is a tree on $\P$ such that $\Q(\T)$ doesn't exist, i.e., $\T$ has no $\Q$-structures. Such trees are called maximal. Then $\emptyset$-iterability says that $\T$ has a branch $b$ such that  $\M^\T_b$ is full, has all the relevant mice projecting to $\d(\T)$.}
\end{enumerate}
\end{definition}

The following is our main theorem of this section. 

\begin{theorem}\label{main theorem 1} Suppose $\K$ is a translatable structure. Then, in $\K$, there is a proper class premouse $\M$ such that $\M$ has a proper class of Woodin cardinals and a strong cardinal.
\end{theorem}

We spend the rest of this section describing $\M$. Just like in Steel's construction, we will introduce an operator in order to compute extenders, much like the way Steel does in \cite{DMATM}. 

For the rest of this section, we fix a translatable structure $\K$ and work entirely inside $\K$. Let $\d<\l$ be the first two Woodin cardinals of $\K$. Let $\P=\W_\omega(\K|\d)$ and $\N=(\mathcal{J}^{\vec{E}}[\P])^{\K|\l}$. Thus, $\N$ is the output of the $\mathcal{J}^{\vec{E}}$ construction of $\K|\l$ done over $\P$ using extenders with critical point $>o(\P)$\footnote{In this paper, we let $o(X)=_{def}\sup(X\cap Ord)$.}.

Before we go on, we record the following useful lemmas about translatable structures. Recall from \cite{ATHM} (see Definition 1.26 of \cite{ATHM}) that if $\S$ is a premouse and $\eta<o(\S)$ then 
\begin{center}$\mathcal{O}_\eta^{\S}=\cup\{ \S || \b : \rho_\omega(\S||\b)\leq \eta$ and $\eta$ is a cutpoint in $\S\}$.\end{center}

\begin{lemma}\label{computing successors correctly} Suppose $\T\in \K$ is a normal tree on $\K$ without fatal drops\footnote{Recall from \cite{ATHM} (see Definition 1.27 of \cite{ATHM}) that we say $\T$ has a fatal drop if for some $\a<lh(\T)$ and $\eta$, the rest of $\T$ after stage $\a$ is based on $\mathcal{O}_\eta^{\M^\T_\a}$ and is above $\eta$.} based on the window $(\d, \l)$ such that $\T$ has a limit length $\leq \l$, $\T$ doesn't have a branch in $\K$ and for every limit $\a<lh(\T)$, $\Q(\T\rest \a)\insegeq \M^\T_\a$ exists\footnote{Recall that $\Q(\T\rest \a)$ is the $\Q$-structure of $\T\rest \a$ (see Definition 1.25 of \cite{ATHM}).} and is $Ord$-iterable as a $\Sigma^\K$-mouse. Then for a stationary many $\a<\l$, $(\a^+)^{\M(\T)}=(\a^+)^\K$. 
\end{lemma}
\begin{proof}
Let $\l^*$ be the sup of the first $\omega$ Woodin cardinals of $\K$ and let $g\subseteq Coll(\omega, <\l^*)$ be $\K$-generic. Let $M$ be the derived model computed in $\K(\bR^*)$ and let $\Phi=\Sigma^g\rest HC^M$. We let $\R$ be the $\Phi$-mouse representation of $(V_\Theta^{\H_\Phi})^M$. More precisely, $\M$ is the direct limit of pairs $(\Q, B)$ such that $\Q$ is a suitable $B$-iterable $\Phi$-premouse and $B\in OD_{\Phi}\cap \powerset(\bR)$. For more on this construction see Section 1.3 of \cite{TrangSarg}. 

It follows from clause 5 of \rdef{translatable structure} that there is an embedding $\pi:Hull_1^{\K|(\l^{+\omega})}\rightarrow Hull_1^\R$. Thus, it is enough to prove that in $\R$, whenever $\T$ is a normal tree on $\R$ such that $lh(\T)\leq \Theta^M$, $\T$ has no fatal drops and $\T$ is correctly guided\footnote{This notion can be found in \cite{ATHM} on page 105. It essentially says that $\T$ chooses unique branches identified by iterable $\Q$-structures.} then either $\T$ has a branch or for stationary many $\a<\Theta^M$, $(\a^+)^{\M(\T)}=(\a^+)^\R$. It follows from Woodin's computation of $\H$ relativized to $\Phi$, that $\H^M_\Phi\models ``\R$ is iterable for trees that are in $\R$" (see, for instance, Theorem 1.3 of \cite{Trang}).

The proof now follows easily. Fix  a normal correctly guided $\T\in \R$ on $\R$ such that $lh(\T)\leq \Theta^M$. Let $b\in \H^M_\Phi$ be the branch of $\T$. We have two cases. Suppose first that $\pi^\T_b[\Theta^M]\not \subseteq \Theta^M$. It then easily follows that $\Q(\T)$ exists and $\Q(\T)\insegeq \M^\T_b$. But because $\Q(\T)\in \R$, we have that $b\in \R$. It follows that we must have that $\pi^\T_b[\Theta^M]\subseteq \Theta^M$. It now follows that the set of $\a<\Theta$ such that $(\a^+)^\R=(\a^+)^{\M(\T)}$ contains a club in $\H^M_\Phi$. Therefore, that set is stationary in $\R$.
\end{proof}

Following Steel's translation procedure, we will describe an operator $S:\K|\l \rightarrow \K|\l$\footnote{$S$ stands for Steel.}. The domain of our operator $S$ will consist of $X\in \K|\l$ such that $X$ is a premouse with the property that $\P\insegeq X$. To define the $S$-operator over one such $X$, we first need to isolate a large class of \textit{universal models} over $X$. 

\subsection{Universal mouse} 

Fix a premouse $X\in \K|\l$ such that $\P\insegeq X$. The universal mice over $X$ that we will consider act like the mice that are constructed via fully backgrounded constructions in $\K|\l$. Suppose $\eta<\l$ is such that $X\in \K|\eta$. Let then $\Q^{X, \eta}$ be the output of $\mathcal{J}^{\vec{E}}[X]$-construction of $\K|\l$ in which all extenders used have critical points $>\eta$. we have that $L[\Q^{X, \eta}]\models ``\l$ is a Woodin cardinal".  It follows from clause 3 of \rdef{translatable structure} that 
\begin{lemma}
 $\Q^{X, \eta}$ is $(\l, \l)$-iterable. 
 \end{lemma}
 We let $\Lambda^{X, \eta}$ be the $(\l, \l)$-strategy induced by the background strategy. 

\begin{definition}[Universal mouse]\label{universal model}
 We say $\Q\subseteq \K|\l$ is a universal mouse over $X$ if $\Q$ is $(\l, \l)$-iterable, for every $\eta\in (o(X), \l)$, $\mathcal{O}^\Q_\eta\insegeq \W(\Q|\eta)$ and for a stationary set of $\xi<\l$, $(\xi^+)^\Q=\xi^+$. 
\end{definition}

Our first lemma on universal models shows that they compute lower parts correctly.
 
\begin{lemma}\label{lp persistence} 
Suppose $\Q$ is a universal mouse over $X$. Then the function $Y\rightarrow \W(Y)$ is definable over $\Q$.
\end{lemma}
\begin{proof}
We first show that if $\eta>o(X)$ is a $\Q$-cardinal then $\mathcal{O}^\Q_\eta=\W(\Q|\eta)$. To see this, let $E\in \vec{E}^\Q$ be the least extender such that $\eta$ is a cutpoint in $Ult(\Q, E)$. It is now enough to show that $Ult(\Q, E)|(\eta^+)^{Ult(\Q, E)}=\W(\Q|\eta)$. Since we already have that 
\begin{center}
$Ult(\Q, E)|(\eta^+)^{Ult(\Q, E)}\insegeq \W(\Q|\eta)$,
\end{center}
 it is enough to show that $\W(\Q|\eta)\insegeq (\eta^+)^{Ult(\Q, E)}$. 

Fix then an $\M\insegeq \W(\Q|\eta)$ such that $\rho_\omega(\M)=\eta$ and $\M\not \insegeq Ult(\Q, E)$. Let $\pi: H\rightarrow \K|(\l^+)^\K$ be such that letting $\cp(\pi)=\nu$, $\nu$ is an inaccessible cardinal in $\K$, $\pi(\nu)=\l$, $E\in \vec{E}^{\Q|\nu}$ and $Ult(\Q, E)\in rng(\pi)$\footnote{This reflection into $H$ is required because of limited iterability.}. Let $\Lambda$ be the $(\l, \l)$-strategy of $Ult(\Q, E)$ and let $\T$ and $\U$ be the trees of length $\nu+1$ on respectively $Ult(\Q, E)|\nu$ and $\M$ build using the comparison process. Let $\R$ be the last model of $\T$ and $\S$ be the last model of $\U$. 

Because $\M$ must win the comparison with $Ult(\Q, E)|\nu$, we have that $\pi^\T[\nu]\subseteq \nu$. It follows that $\R$ computes unboundedly many successors below $\nu$ correctly, implying that $\M$ couldn't out iterate $\R$, contradiction!

Next, fix $Y\in \Q$. We want to see that $\W(Y)$ is definable. Let $\eta$ be a cardinal such that $Y, X\in \Q|\eta$. \\

\textit{Claim.}
$\M\insegeq\W(Y)$ if an only if there is $\R\insegeq \mathcal{O}^\Q_\eta$ such that $\rho(\R)=\eta$, $\R$ has $\omega$ many Woodin cardinals with sup $\nu$, $\M\in \R$ and $\R\models  ``\M$ is $\omega_1$-iterable in the derived model at $\nu$".\\\\
\begin{proof}
We prove the forward direction. The reverse direction follows from the homogeneity of the collapse. Let $\l^*$ be the sup of the first $\omega$ many Woodin cardinals of $\K$ and let $g\subseteq Coll(\omega, <\l^*)$ be $\K$ generic. Let $M$ be the derived model computed in $\K(\bR^*)$. It follows from clause 2 of \rdef{translatable structure} that $\mathcal{O}^\Q_\eta=Lp(\Q|\eta)$ (see \cite[Theorem 3.10]{ATHM}). We also have that $\M\insegeq Lp(Y)$. 

Let $\Phi=\Sigma^g\rest HC^M$. Because $\Phi\in M$, we have that $M\models \theta_0<\Theta$. It follows that in $M$, there is a pair $(\S, \Lambda)$ such that $\S$ is $\omega$-suitable mouse over $\Q|\eta$, $\Lambda$ is a fullness preserving iteration strategy for $\S$ and every $(\mathbf{\Sigma}^2_1)^M$ set is in the derived model of $\S$ as computed via $\Lambda$ (see, for instance, the proof of Lemma 5.13 of \cite{ATHM}). 

We then have that $\M\in \S$ and letting $\xi$ be the sup of the Woodin cardinals of $\S$, $\S\models ``\M$ has an $\omega_1$-iteration strategy in the derived model at $\xi$". We can then obtain such an $\R$ by simply taking a hull of $\S$ which is $\eta$-sound. 
\end{proof}

Clearly the claim finishes the proof of the lemma. 
\end{proof}

\subsection{Good universal mouse}

We continue working with our $X\in \K|\l$. Given any $X$-premouse $\Q\subseteq \K|\l$ we let $\k^\Q$ be the least $<\l$-strong of $\Q$ if exists. Suppose $\eta\in (o(X), \l)$. 

\begin{lemma}
$\k^{\Q^{X, \eta}}$ is a limit of Woodin cardinals in $\Q^{X, \eta}$.
\end{lemma}
\begin{proof}
Let $\Q=\Q^{X, \eta}$. Assume then that for some $\eta$, $\Q=\Q^{X, \eta}$. We show that $\l$ is a limit of $\Q$-cardinals $\xi$ such that $\W(\Q|\xi)\models ``\xi$ is a Woodin cardinal". The proof will in fact show that there is a club of such $\xi$. Suppose not and let $\eta$ be the sup of such $\xi$. Let $\S$ be the output of $\mathcal{J}^{\vec{E}}$ construction of $\Q$ in which all extenders used have critical point $>\eta$. Then standard arguments show that $\P$ iterates to $\S$ (see Lemma 2.11 of \cite{ATHM})\footnote{If $\P$ side comes short then we can use $\S$-constructions to show that there is $\xi>\eta$ with the property that $\W(\Q|\xi)\models ``\xi$ is a Woodin cardinal".}. But this implies that $\S$ is not universal contradicting clause 4 of \rdef{translatable structure} (see Lemma 2.13 of \cite{ATHM}).
\end{proof}

Fix now $\eta\in (o(X), \l)$ and let $\Q=\Q^{X, \eta}$. Just like in the proof of the \textit{Mouse Set Conjecture}, we would like to recover a tail of $(\P, \Sigma)$ inside $\Q$. The details of this construction can be found in Section 6.3 of \cite{ATHM}. It has also appeared in Section 4 of \cite{TrangSarg}. Below we just review the construction. 

Suppose $g\subseteq Coll(\omega, <\k^\Q)$ is a $\Q$-generic and $M$ is the derived model computed in $\Q[g]$.
Working in $M$, let $\mathcal{F}=\{ (\R, A):\R$ is suitable and $A$-iterable$\}$. Let $\M_\infty^\Q$ be the direct limit of $\mathcal{F}$. It then follows that $\M_\infty^\Q\in\Q$ and letting $\xi$ be the Woodin cardinal of $\M_\infty^\Q$, 
\begin{center}
$\W_\omega(\M_\infty^\Q|\xi)=\M_\infty^\Q$.
\end{center}
We let
\begin{center}
 $\R^\Q=\M_\infty^\Q|(\xi^{+\omega})^{\M_\infty^\Q}$.
 \end{center}
Notice that we could define $\R^\S$ for any $\S$ that has the same properties as $\Q$. In particular it is defined for all iterates of $\Q$. The following lemma can be proved using the proof of Theorem 6.5 of \cite{ATHM} (see the successor case).

\begin{lemma}\label{msc proof}
 Let $\Lambda$ be the $(\l, \l)$-iteration strategy for $\Q$ induced by the background strategy. Let $i:\Q\rightarrow \M$ be an iteration according to $\Lambda$. Then
$\R^\M$ is a $\Sigma$-iterate of $\P$ and whenever $g$ is a $\K$-generic for a poset in $\M$, 
\begin{center}
$\Sigma^{g}_{\R^\M}\rest \M[g]\in \mathcal{J}[\M[g]]$. 
\end{center}
\end{lemma}

The next lemma shows that all iterates of $\Q^{X, \eta}$ are universal.

\begin{lemma}\label{background models are universal} Suppose $\eta\in (o(X), \l)$ and that $\Lambda$ is the strategy of $\Q^{X, \eta}$ induced by the $(\l, \l)$-strategy of $\K|\l$. Let $\M$ be a $\Lambda$-iterate of $\Q^{X, \eta}$ such that the iteration embedding $i:\Q^{X, \eta}\rightarrow \M$ exists. Then $\M$ is a universal mouse over $X$.
\end{lemma}
\begin{proof}
It is enough to show that for stationary many $\a$, $\M$ computes $(\a^+)^\K$ correctly (as this implies that $\M$ wins coiteration with any other mouse). Let $\K^*$ be the iterate of $\K$ obtained by applying $\P$-to-$\R^\M$ iteration to $\K$. Let $\K^{**}=(\mathcal{J}^{\vec{E}, \Sigma_{\R^\M}})^\M$. It follows from universality of background constructions (see Lemma 2.13 of \cite{ATHM}) that $\K^*$ and $\K^{**}$ co-iterate to a common $\K^{***} \subseteq \K|\l$ via normal trees $\T$ and $\U$ respectively. 

Suppose first $\T$ has a branch $b\in \K$. It follows from the universality of $\K^{**}$ that $\pi^\T_\b[\l]\subseteq \l$. This then implies that for stationary many $\a$, $\K^{***}$, and hence $\K^{**}$,  compute $(\a^+)^\K$-correctly. If $\T$ doesn't have a branch then it follows from \rlem{computing successors correctly} that for stationary many $\a$, $\K^{***}$ and hence $\K^{**}$,  compute $(\a^+)^\K$-correctly.
\end{proof}

We now define the \textit{good universal mice}.

\begin{definition}[Good universal mouse]\label{good universal models} We say $\Q$ is a good universal mouse if it has a $(\l, \l)$-iteration strategy $\Lambda$ such that whenever $i:\Q\rightarrow \M$ is an iteration according to $\Lambda$, $\M$ is universal,  
$\R^\M$ is a $\Sigma$-iterate of $\P$ and whenever $g$ is a $\K$-generic for a poset in $\M$, 
\begin{center}
$\Sigma^{g}_{\R^\M}\rest \M[g]\in \mathcal{J}[\M[g]]$. 
\end{center}

\end{definition}

\rlem{msc proof} and \rlem{background models are universal} imply that there are good universal mice. Let then $\Q$ be a good  universal mouse. We then let $\pi^\Q: \P\rightarrow \R^\Q$ be the iteration embedding according to $\Sigma$. 

 From now on we let $A\subseteq \P$ be the collection of all terms that $\Sigma$ guides correctly. For $\tau\in A$ we let $\pi^\Q_\tau=\pi^\Q\rest H_\tau^\P$ and set
\begin{center}
$\Q^+=(\Q|((\k^\Q)^+)^\Q, \in, \{\pi^\Q_\tau\}_ {\tau\in A})$.
\end{center} 
Finally we can define the $S$-operator.
\begin{definition}\label{s operator} Suppose $X\in \K|\l$ and $\Q$ is a good universal mouse over $X$. Then
\begin{center}
$S^\Q(X)=Hull_1^{\Q^+}(X)$.
\end{center}
\end{definition}

In order to prove useful lemmas about the $S$ operator, we use a standard way of splitting the $S$ operator into smaller pieces. 
For a finite $B \subseteq A$, we set
\begin{center}
$\Q^+_B=(\Q|((\k^\Q)^+)^\Q, \in, \{\pi^\Q_\tau: \tau\in B\} )$.
\end{center}
Let then
\begin{center}
$S_B(X)=Hull^{\Q^+_B}_1(X)$ 
\end{center}
Clearly if $B\subseteq C\subseteq A$ are such that $B$ and $C$ are finite then there is a canonical $\sigma^{X}_{B, C}: S_B(X)\rightarrow S_C(X)$. It then follows that the direct limit of $(S_B: B\subseteq_{fin}A)$ under the maps $\sigma^X_{B, C}$ is just $S(X)$.

Our first lemma shows that the $S$ operator is independent of the choice of $\Q$. 

\begin{lemma} $S(X)$ is independent of the choice of $\Q$. 
\end{lemma}
\begin{proof}
Let $\Q_0$ and $\Q_1$ be two good universal mice over $X$. Let $\pi: H\rightarrow \K|(\l^+)^\K$ be such that if $\cp(\pi)=\nu$ then $\nu$ is an inaccessible cardinal in $\K$, $X\in H$, $\pi(\nu)=\l$ and $\Q_0, \Q_1\in rng(\pi)$. Moreover, letting for $i=0,1$, $\bar{\Q}_i=\pi^{-1}(\Q_i)$, for each $i\in 2$ there are unboundedly in $\nu$ cardinals $\xi$ such that 
\begin{center}
$(\xi^+)^{\bar{\Q}_i}=(\xi^+)^\K$\ \ \ \ \ \ \ \ (*).
\end{center}
 Notice that we have that for $i=0,1$, $\bar{\Q}_i$ is fully iterable (see clause 3 of \rdef{translatable structure}). Let $\Lambda_i$ be an iteration strategy for $\bar{\Q}_i$. Comparing $\bar{\Q}_0$ with $\bar{\Q}_1$ using $\Lambda_0$ and $\Lambda_1$ respectively, we obtain $(\U_0, \Q_3)$ and $(\U_1, \Q_4)$ such that
\begin{enumerate}
\item $\U_i$ is a tree on $\bar{\Q}_i$ according to $\Lambda_i$,
\item $\Q_3$ is the last model of $\U_0$,
\item $\Q_4$ is the last model of $\U_1$,
\item $\Q_3\insegeq \Q_4$ or $\Q_4\insegeq \Q_3$.
\end{enumerate}
It follows from the usual comparison argument and (*) that $\Q_3=\Q_4$ and $o(\Q_3)=\nu$. 

Suppose now $\tau\in A$. Then
\begin{center}
 $\pi^{\U_0}[\pi_\tau^{\Q_0}]=\pi^{\U_1}[\pi_\tau^{\Q_1}]=_{def}\sigma_\tau$.
 \end{center}
It then easily follows that if $\xi$ is the least $<\nu$-strong of $\Q_3$ then letting $D=(\Q_3|(\xi^+)^{\Q_3}, \in, \{\sigma_\tau: \tau\in A\})$,
\begin{center}
$S^{\Q_0}(X)=\S^{\Q_1}(X)=Hull^D_1(X)$.
\end{center}

\end{proof}

\subsection{Basic properties of $S(X)$}

The following is an easy lemma and follows from super fullness preservation of $\Sigma$.
\begin{lemma}
$\W(X)\in \S(X)$
\end{lemma}

Next we would like to show that $S$ operator is somehow easy to compute relative to $\Sigma$. Our first theorem towards it is that $\Sigma$ moves $S$ operator correctly. For the duration of this section let $\nu$ be the second Woodin cardinal of $\N$ and let $\Lambda$ be the iteration strategy of $\N$ that acts on stacks that based on the window $(\d, \nu)$ (recall that $\N=(\mathcal{J}^{\vec{E}}[\P])^{\K|\l}$).

\begin{lemma}
Suppose $i:\N\rightarrow \bar{\N}$ is an embedding coming from an iteration according to $\Lambda$. Then for any $B\subseteq_{fin} A$,
\begin{center}
$i(S_B(\N|\nu))=S_B(\bar{\N}|i(\nu))$
\end{center}
\end{lemma}
\begin{proof}
The claim follows by noticing that $\N$ is a universal mouse over $\N|\nu$ and $\bar{\N}$ is a universal mouse over $\bar{\N}|i(\nu)$. It also uses the following claim.\\

\textit{Claim.} For $\tau\in A$, $i(\pi^\N_\tau)=\pi^{\bar{\N}}_\tau$.\\\\
\begin{proof}
We have that there is $j:\K\rightarrow \K^*$ and $\sigma:\bar{\N}\rightarrow j(\N)$ such that $j\rest \N=\sigma\circ i$. We then have that for every $\tau\in A$, $j(\pi^\N_\tau)=\pi^{j(\N)}_\tau$. The claim now follows from clause 2 of \rdef{translatable structure}.
\end{proof}
\end{proof}

Our next lemma shows that models obtained via $\S$-constructions from $S$ operators are themselves $S$ operators.

\begin{lemma}\label{transferring s-operator} Suppose $i:\N\rightarrow \bar{\N}$ is an embedding coming from an iteration according to $\Lambda$. Suppose $X$ is such that $X$ is generic over $\mathcal{J}_\omega[\bar{\N}|i(\nu)]$ and $\mathcal{J}_\omega[\bar{\N}|i(\nu)][X]=\mathcal{J}[X]$. Fix $B\subseteq_{fin}  A$ and let $D$ be the mouse obtained over $X$ from $S_B(\bar{\N}|i(\nu))[X]$ via $\S$-construction. Then 
\begin{center}
$D=S_B(X)$.
\end{center}
\end{lemma}
\begin{proof}
The proof is a straightforward chase thru definitions. The important fact to note is that if $\N^*$ is the $X$-mouse obtained from $\bar{\N}[X]$ via $\S$-construction then $\N^*$ is a good universal mouse over $X$, $\R^{\N^*}=\R^{\bar{\N}}$ and $\pi^{\N^*}=\pi^{\bar{\N}}$.
\end{proof}

\subsection{Certified extenders and premice}

Next we show that the $S$-operator is captured by $\M_1^{\#, \Sigma}$ operator. Recall that $\nu$ is the second Woodin cardinal of $\N$. Let $\k=\k^\N$.

\begin{lemma}\label{towards uniformity} Suppose $\xi\in ((\kappa^+)^\N, \l)$ is an $\N$-cardinal. Then $S(\N|\xi)\in \M_1^{\#, \Sigma}(\N|\xi)$.
\end{lemma}
\begin{proof}
 Let $\N^*=Ult(\N, E)$ where $E\in \vec{E}^\N$ is such that $Ult(\N, E)\models ``\xi$ is a cutpoint". Notice that $\N^*$ is a universal mouse over $\N|\xi$. Let $i:\N\rightarrow \bar{\N}$ be iteration of $\N$ according to $\Lambda$ such that $\N|\xi$ is generic over the extender algebra of $\bar{\N}$ at $i(\nu)$. Let $\Q=\bar{\N}|(i(\nu)^{+\omega})^{\bar{\N}}$. Notice that $\N^*|(\xi^{+\omega})^{\N^*}$ is generic over $\mathcal{J}_{\omega}[\Q]$ and 
 \begin{center}
 $\mathcal{J}_{\omega}[\Q][\N^*|(\xi^{+\omega})^{\N^*}]=\mathcal{J}_{\omega}[\N^*|(\xi^{+\omega})^{\N^*}]$.
 \end{center}
 It follows from the proof of \rlem{transferring s-operator} that if $D$ is the model obtained via $\S$-construction from $S(\Q)[\N^*|(\xi^{+\omega})^{\N^*}]$ then $D=S(\N^*|(\xi^{+\omega})^{\N^*})$. Notice also that
 \begin{center}
 $S(\N|\xi)=Hull^{S(\N^*|(\xi^{+\omega})^{\N^*})}(\N|\xi)$.
 \end{center}
 
 Finally, notice that $\Q\in \M_1^{\#, \Sigma}(\N|\xi)$, $S(\N|(\nu^{+\omega})^\N)\in \M_1^{\#, \Sigma}(\N|\xi)$ and 
 \begin{center}
 $S(\Q)=\cup_{B\subseteq_{fin} A}i(S_B(\N|(\nu^{+\omega})))$.
 \end{center}
  It then follows that $S(\N|\xi)\in \M_1^{\#, \Sigma}(\N|\xi)$.
\end{proof}

The next lemma shows that extender can be uniformly identified via the $S$-operator.

\begin{lemma}\label{uniformity}
Suppose $\xi\in [(\k^+)^\N, \l)$ be such that for some $E\in \vec{E}^\N$ with $\cp(E)=\k$, $E$ is indexed at $\xi$. Let $\zeta$ be least cutpoint $>\k$ of $Ult(\N, E)$. Then $E\in \M_1^{\#, \Sigma}(Ult(\N, E)|\zeta)$ and has a uniform definition from $Ult(\N, E)|\zeta$ and $\Sigma$.
\end{lemma}
\begin{proof}
Notice that
\begin{center}
$S(\N|\k)=Hull_1^{S(Ult(\N, E)|\zeta)}(\N|\k)$\\
$S(Ult(\N, E)|\zeta)=Hull_1^{\pi_E(S(\N|\k))}(Ult(\N, E)|\zeta)$
\end{center}
Let then $\sigma: S(\N|\k)\rightarrow S(Ult(\N, E)|\zeta)$ and $k:S(Ult(\N, E)|\zeta)\rightarrow \pi_E(S(\N|\k))$ be the uncollapse maps. We then have that $\pi_E\rest S(\N|\k) =k\circ \sigma$. It follows from \rlem{towards uniformity} that $\sigma\in \M_1^{\#, \Sigma}(Ult(\N, E)|\zeta)$.

It follows from the construction that $B\in E_a$ if and only if $a\in \sigma(B)$. It then also follows that $E\in \M_1^{\#, \Sigma}(Ult(\N, E)|\zeta)$. Chasing thru the proof of \rlem{towards uniformity} and the current proof, we easily see that the definition of $E$ in $\M_1^{\#, \Sigma}(Ult(\N, E)|\zeta)$ is uniform from $Ult(\N, E)|\zeta$ and $\Sigma$.
\end{proof}

We now give a more precise description of the procedure that defines the extender on the sequence of $\N$. Recall that we let $\k=\k^\N$. 

\begin{definition}[$\Sigma$-certified extenders]\label{certified extenders} Suppose $\Q$ is a premouse such that $\N|\k\insegeq \Q$ and $\N|(\k^+)^\N=\Q|(\k^+)^\Q$. We say $E\in \M_1^{\#, \Sigma}$ is a $(\Sigma, \Q)$-certified extender if 
\begin{enumerate}
\item $\cp(E)=\k$, 
\item $\Q\inseg\pi_E^\Q(\Q|\k)$,
\item $Hull_1^{S(\Q)}(\N|\k)=S(\N|\k)$,
\item $E$ is the $(\k, o(\Q))$-extender derived from $k:S(\N|\k)\rightarrow S(\Q)$ where $k$ is the uncollapse map. 
\end{enumerate}
We say $\Q$ is a $\Sigma$-certified premouse if $\N|\k\insegeq \Q$, $\N|(\k^+)^\N=\Q|(\k^+)^\Q$ and whenever $E_\a\in \vec{E}^\Q$ is an extender such that $\cp(E_\a)=\k$, $E_\a$ is $(\Sigma, \Q|\a)$-certifed.
\end{definition}

The next lemma is a condensation lemma.

\begin{lemma}\label{condensation} Suppose $\Q$ is a $\Sigma$-certified premouse, and $\R$ is such that $\N|\k\insegeq \R$ and $\R|(\k^+)^\R=\N|(\k^+)^\N$. Suppose that there is an embedding $i:\R\rightarrow \Q$. Then $\R$ is $\Sigma$-certified.
\end{lemma}
\begin{proof}
Suppose $\a\in dom(\vec{E})$ is such that $\R|\a$ is $\Sigma$-certified but $E_\a^\R$ is not $(\Sigma, \R|\a)$-certified. Let $E_\b^\Q=i(E_\a)$. Let $\sigma:Ult(\R, E_\a^\R)\rightarrow Ult(\Q, E_\b^\Q )$ be the canonical factor map. Let $\M_0\in Ult(\Q, E_\b^\Q)$ and $\M_1\in Ult(\R, E_\a^\R)$ be the iterates of $\W_\omega(\N|\nu)$ (recall that $\nu$ is the second Woodin cardinal of $\N$) according to $\Lambda$ that make respectively $\Q|\b$ and $\R|\a$ generic. It follows from branch condensation of $\Lambda$ (the strategy of $\N$ based on the window $(\d, \nu)$) and clause 2 of \rdef{translatable structure} that if for $i\in 2$, $\pi_i:\W_\omega(\N|\nu)\rightarrow \M_i$ are the iteration embeddings according to $\Lambda$, then $\pi_0=(\sigma\rest \M_1)\circ \pi_1$. It then follows from the proof of \rlem{towards uniformity} that $\sigma\rest \R|\a$ can be lifted to $\sigma^+: S(\R|\a)\rightarrow S(\Q|\b)$ such that letting $k:S(\N|\k)\rightarrow S(\Q|\b)$ and $l: S(\N|\k)\rightarrow S(\R|\a)$ be the uncollapse maps,
\begin{center}
$k=\sigma^+\circ l$.
\end{center}
But then we have that 
\begin{eqnarray*}
(a, B)\in E_\a^\R &\iff&
a\in \pi_{E_\a^\R}(B)\\ &\iff&  \sigma(a)\in \pi_{E_\b^\Q}(B)\\ &\iff& \sigma(a)\in k(B)\\ &\iff& \sigma^+(a)\in \sigma^+(l(B))\\ &\iff& a\in l(B)
\end{eqnarray*}
\end{proof}

\subsection{The construction}

In the section we describe a construction that extends $\N=(\mathcal{J}^{\vec{E}})^{\K|\l}$ and converges to the model $\M$ as in \rthm{main theorem 1}. The construction is a kind of $K^c$ construction extending $\N$. In this construction, there are two conditions on the background extenders. Suppose $E$ is an extender that we would like to use as a background extender. Suppose first that $\cp(E)>\l$. In this case, we require that $E\in \vec{E}^{\K}$ and satisfy the background conditions imposed by the usual fully backgrounded $\mathcal{J}^{\vec{E}}$-constructions. Suppose next that $\cp(E)<\l$. In this case, it is required that $\cp(E)=\k^\N$ and moreover, letting $\Q$ be the stage of the construction where it is possible to add $E$, it is required that $E$ is $(\Sigma, \Q)$-certified. The following is a more precise description of the construction.

\begin{definition}\label{the construction of the model}
Let $(\M_\xi, \N_\xi : \xi<\Omega)$ be a sequence of premice defined as follows:
\begin{enumerate}
\item $\N_0=\N$.
\item $\M_\xi=(\mathcal{J}_\a^{\vec{E}}, \in, \vec{E})$ is a passive premouse. $\N_{\xi+1}$ is obtained from $\M_\xi$ as follows:
\begin{enumerate}
\item There is an extender $F^*$ on the sequence of $\K$ with critical point $>\l$ such that $F^*$ coheres the construction and for some $\nu<\a$, if $F\rest \nu=F^*\cap ( [\nu]^{<\omega}\times \mathcal{J}_\a^{\vec{F}})$ then
\begin{center}
$\N_{\xi+1}=(\mathcal{J}_\a^{\vec{E}}, \in, \vec{E}, E)$
\end{center}
\item The condition above fails, and there is a $(\Sigma, \M_\xi)$-cerified extender $E$ such that  
\begin{center}
$\N_{\xi+1}=(\mathcal{J}_\a^{\vec{E}}, \in, \vec{E}, E)$
\end{center}
is a premouse.
\item Both (a) and (b) fail and $\N_{\xi+1}=\mathcal{J}(\M_\xi)$.
\end{enumerate}
In all the above cases, we let $\M_{\xi+1}=\mathbb{C}_\omega(\N_{\xi+1})$.
\item Suppose $\l$ is limit. Then $E=E_\a^{\M_\l}$ for some $\a$ if and only if there is $\xi<\l$ such that for all $\b\in (\xi, \l)$, $E=E_\a^{\M_\xi}$.
\end{enumerate}
The construction halts at stage $\xi\in Ord$ if either $\mathcal{C}_\omega(\N_\xi)$ is not defined or else $\rho_\omega(\M_\xi)<\l$. We say the construction converges if doesn't halt at any $\xi$. 
\end{definition}

Our next goal is to show that the construction converges. A similar argument then will show that $\k^\N$ is a strong cardinal in the final model. Proofs of both of these facts are very similar to the proofs given by Steel in \cite{DMATM} (see Section 13 and 15 of that paper). Our first lemma shows that certain comparisons do not encounter disagreements with extenders that have critical point $\k^\N$.

\begin{lemma}\label{no bad disagreement} 
Suppose $\l<\xi<\theta$ are two ordinals. Suppose further that $\pi: H\rightarrow \K|\theta$ is such that letting $\cp(\pi)=\nu$, the following holds:
\begin{enumerate}
\item $(\M_\a, \N_\a: \a<\xi) \in rng(\pi)$ and $\N_\xi\in rng(\pi)$.
\item $\nu\in (\k^\N, \l)$ is an inaccessible cardinal in $\K$, and $\pi(\nu)=\l$.
\item For every $\a\leq \xi$, $\N_\a\models ``\l$ is a Woodin cardinal".
\end{enumerate}
Let $\bar{\N}=\pi^{-1}(\N)$, $(\S_\a, \Q_\a: \a<\zeta)=\pi^{-1}((\M_\a, \N_\a: \a<\xi))$ and $\Q_\zeta=\pi^{-1}(\N_\xi)$. Let $\S\insegeq \N$ be the least such that $\S\models ``\nu$ is a Woodin cardinal" and $\mathcal{J}_1(\S)\models ``\nu$ is not a Woodin cardinal". Then the comparison of the construction of $H$ that produces $(\S_\a, \Q_\a: \a<\zeta)$ and $\Q_\zeta$ with $\S$ doesn't encounter disagreements involving extenders with critical point $\k^\N$.
\end{lemma}
\begin{proof}
Suppose not. Let $\T$ and $\U$ be the iteration trees on $H$ and $\S$ respectively coming from the comparison process that expose a disagreement involving extenders with critical point $\k^\N$. Let $H^*$ and $\S^*$ be the last models of $\T$ and $\U$ respectively. Notice that both $\T$ and $\U$ are above $\nu$. 

Let $\Q$ be the model appearing in the construction of $H^*$ before stage $\pi^\T(\zeta)$ such that the disagreement between the construction of $H^*$ and $\S^*$ is between $\Q$ and $\S^*$. Suppose first that there is an extender with critical point $\k^\N$ on the sequence of $\Q$ that causes a disagreement. Let $\a$ be such that $E_\a^\Q$ is this extender. Thus we have that $\Q|\a=\S^*|\a$.

Suppose $E_\a^{\S^*}\not=\emptyset$. It then follows from our construction, \rlem{condensation} and \rlem{uniformity} that both $E_\a^{\S^*}$ and $E_\a^\Q$ are $(\Sigma, \Q|\a)$-certified. Clause 4 of \rdef{certified extenders} implies that $E_\a^{\S^*}=E_\a^\Q$. 

Suppose next that $E_\a^{\S^*}=\emptyset$. We then have two cases. Suppose first for some $\b>\a$, $E_\b^{\S^*}\not=\emptyset$ and $\cp(E_\b^{\S^*})=\k$. Fix minimal such $\b$. Then $\a$ is a cutpoint cardinal of $Ult(\S^*, E_\b^{\S^*})$ and $(\k, \a)$-extender derived from $\pi_{E_\b^{\S^*}}$ is $(\Sigma, \Q|\a)$-certified. Then again Clause 4 of \rdef{certified extenders} and initial segment condition imply that $E_\a^{\S^*}=E_\a^\Q$. 

Next suppose that $\a$ is a cutpoint of $\S^*$. Because $E_\a^\Q$ is $(\Sigma, \Q|\a)$-certified, we have that $\W(\Q|\a)\insegeq Ult(\Q, E_\a^\Q)$. It then follows that $\S^*\insegeq Ult(\Q, E_\a^\Q)$ implying that $\Q\models ``\nu$ is not a Woodin cardinal" of $\Q$, contradiction.

Next, suppose that $H^*$ side is not active. Let $\a$ be least such that $\S^*|\a=\Q|\a$, $E_\a^{\S^*}\not =\vec{E}^\Q$ and $\cp(E_\a^{\S^*})=\k^\N$. Arguing as above, we see that $\a$ must be a cutpoint in $\Q$. It then follows that $\Q\insegeq Ult(\S^*, E_\a^{\S^*})$, implying that $\Q\insegeq \S^*$. This contradicts our assumption that $\Q$ participates in a disagreement. 
\end{proof}

We now prove that 

\begin{lemma}\label{defined} For each $\a$, $\M_\a$ is defined. 
\end{lemma}
\begin{proof} Suppose not. Let $\a$ be least such that $\M_{\a}$ is not defined. \\

\textit{Claim.} $\N_\a\models ``\l$ is a Woodin cardinal".\\\\
\begin{proof}
Suppose not.  Let $\xi$ be least such that $\N_\xi\models ``\l$ is a Woodin cardinal" but $\mathcal{J}_1(\N_\xi)\models ``\l$ isn't a Woodin cardinal". Fix $\theta$ such that the construction up to stage $\xi$ can be done inside $\K|\theta$. Let $\pi: H\rightarrow \K|\theta$ be such that letting $\cp(\pi)=\nu$, the following holds:
\begin{enumerate}
\item $(\M_\a, \N_\a: \a<\xi) \in rng(\pi)$ and $\N_\xi\in rng(\pi)$.
\item $\nu\in (\k^\N, \l)$ is an inaccessible cardinal in $\K$, and $\pi(\nu)=\l$.
\end{enumerate}
Let $\bar{\N}=\pi^{-1}(\N)$, $(\S_\a, \Q_\a: \a<\zeta)=\pi^{-1}((\M_\a, \N_\a: \a<\xi))$ and $\Q_\zeta=\pi^{-1}(\N_\xi)$. 
Let $\S\insegeq \N$ be the least such that $\S\models ``\nu$ is a Woodin cardinal" and $\mathcal{J}_1(\S)\models ``\nu$ is not a Woodin cardinal". It follows from \rlem{no bad disagreement} that $\Q_\zeta=\S$. 

Since $\nu$ is an inaccessible cardinal of $\K$, it follows that $\rho_\omega(\S)=\nu$. Since $\mathcal{J}_1(\N_\xi)\models ``\l$ is not a Woodin cardinal", we must have that $H\models ``\nu$ is not a Woodin cardinal", contradiction!
\end{proof}

Notice now that $\a$ must be a successor ordinal. Fix $\theta$ such that the construction up to stage $\a$ can be done inside $\K|\theta$. Let $\pi: H\rightarrow \K|\theta$ be such that letting $\cp(\pi)=\nu$, the following holds:
\begin{enumerate}
\item $(\M_\b, \N_\b: \b<\a) \in rng(\pi)$ and $\N_\a\in rng(\pi)$.
\item $\nu\in (\k^\N, \l)$ is an inaccessible cardinal in $\K$, and $\pi(\nu)=\l$.
\end{enumerate}
Let $\bar{\N}=\pi^{-1}(\N)$, $(\S_\b, \Q_\b: \b<\zeta)=\pi^{-1}((\M_\b, \N_\b: \b<\a))$ and $\Q_\zeta=\pi^{-1}(\N_\a)$. 
Let $\S\insegeq \N$ be the least such that $\S\models ``\nu$ is a Woodin cardinal" and $\mathcal{J}_1(\S)\models ``\nu$ is not a Woodin cardinal". 

It follows from the Claim and \rlem{no bad disagreement} that the comparison process that compares the construction of $H$ producing $(\S_\b, \Q_\b: \b<\zeta)$ and $\Q_\zeta$ with $\S$ must halt without ever using extenders with critical point $\k^\N$. Since $H$-side must lose, we have that $\mathcal{C}_\omega(\Q_\zeta)$ must be defined, contradiction. 
\end{proof}

Finally we show that $\k^\N$ is a strong cardinal in the model produced via the construction of \rdef{the construction of the model}. The proof is very much like the proof of \rlem{defined}.

\begin{lemma}\label{strong} Suppose $\M$ is the model to which the construction of \rdef{the construction of the model} converges. Then $\M\models ``\k^\N$ is a strong cardinal".
\end{lemma}
\begin{proof}
Suppose not. Let $\eta>\l$ be a cutpoint of $\M$. Let $\xi$ be such that $\N_\xi|(\eta^+)^\M=\M|(\eta^+)^\M$. Fix $\theta$ such that the construction up to stage $\xi$ can be done inside $\K|\theta$. Let $\pi: H\rightarrow \K|\theta$ be such that letting $\cp(\pi)=\nu$, the following holds:
\begin{enumerate}
\item $(\M_\a, \N_\a: \a<\xi) \in rng(\pi)$ and $\N_\xi\in rng(\pi)$.
\item $\nu\in (\k^\N, \l)$ is an inaccessible cardinal in $\K$, and $\pi(\nu)=\l$.
\item $\eta\in rng(\pi)$.
\end{enumerate}
Let $\bar{\N}=\pi^{-1}(\N)$, $(\S_\a, \Q_\a: \a<\zeta)=\pi^{-1}((\M_\a, \N_\a: \a<\xi))$ and $\Q_\zeta=\pi^{-1}(\N_\xi)$. 
Let $\S\insegeq \N$ be the least such that $\S\models ``\nu$ is a Woodin cardinal" and $\mathcal{J}_1(\S)\models ``\nu$ is not a Woodin cardinal". 

We now compare the construction of $H$ producing the sequence $(\S_\a, \Q_\a: \a<\zeta)$ and $\Q_\zeta$ with $\S$. Let $i:H\rightarrow H^*$ and $j:\S\rightarrow \S^*$ be result of this comparison. It follows from \rlem{no bad disagreement} that $\cp(i), \cp(j)>\nu$. We must have that $i(\Q_\zeta)\insegeq \S^{**}$. Let $\gg=i(\pi^{-1}(\eta))$. 

Suppose next that $\gg$ is a cutpoint of $\S^{*}$. It then follows that $\S^{*}\insegeq \W(i(\Q_\zeta)|\gg)$. Since $\M_1^{\#, \Sigma}(i(\Q_\zeta))\in H^*$, we have that $\W(i(\Q_\zeta)|\gg)\in H^*$. It then follows that $\S^{*}\in H^*$ implying that $H^*\models ``\nu$ is not a Woodin cardinal", contradiction. 

We must then have that $\gg$ is not a cutpoint in $\S^{*}$. Let $\a$ be least such that $E=_{def}E_\a^{\S^{*}}$ is the first extender overlapping $\gg$ in $\S^*$. It follows from \rlem{uniformity} that $E$ is $(\Sigma, i(\Q_\zeta)|\gg)$-certified. Again, since $\M_1^{\#, \Sigma}(i(\Q_\zeta))\in H^*$, $H^*\models `` E$ is $(\Sigma, i(\Q_\zeta)|\gg)$-certified". Hence, $E$ should have been on the sequence of $i(\Q_\zeta)$, contradiction.
\end{proof}

\section{A question of Wilson}

In a private conversation Trevor Wilson asked the question. 

\begin{question}[Wilson]\label{trev question}
Assume there are proper class of Woodin cardinals. Suppose the class
\begin{center}
$S=\{ \l: \l$ is a limit of Woodin cardinals and the derived model at $\l$ satisfies $AD^++\theta_0<\Theta\}$
\end{center}
is stationary. Is there a transitive model $\M$ satisfying $ZFC$ such that $Ord\subseteq \M$ and $\M$ has a proper class of Woodin cardinals and a strong cardinal?
\end{question}

We show that the answer is yes.

\begin{theorem} Suppose there are proper class of Woodin cardinals. Assume further that the class $S$ above is stationary. Then there is a transitive model $\M$ satisfying $ZFC$ such that $Ord\subseteq \M$ and $\M$ has a proper class of Woodin cardinals and a strong cardinal.
\end{theorem}
\begin{proof} Given $\l\in S$ let $f(\l)$ be the least cardinal $\eta<\lambda$ such that if $g\subseteq Coll(\omega, \eta)$ then the universal $\Pi^2_1$ set of the derived model at $\lambda$ is witnessed to be $<\l$-uB in $V[g]$. Because $S$ is stationary we have that for some stationary class $S^*\subseteq S$, $f\rest S^*$ is a constant function. We assume, without a loss of generality, that $S^*=S$. Let $g\subseteq Coll(\omega, \eta)$ be $V$-generic. 

Suppose first that there is $\l\in S-\eta$ such that the derived model at $\l$ satisfies $\theta_1<\theta$. Then it follows from unpublished result of Steel (but see the proof of Theorem 15.1 of \cite{DMATM}) that there is a countable active mouse $\M$ such that if $\k$ is the critical point of the last extender of $\M$ then $\M\models ``\k$ is a limit of Woodin cardinals and there is a $<\k$-strong cardinal". Iterating the last extender of $\M$ through the ordinals we obtain a class size model as desired. 

Assume then that for all $\l\in S-\eta$ the derived model at $\lambda$ satisfies $AD^++\theta_1=\Theta$. Similarly, using the results of \cite{ATHM}, we can conclude that for every $\l\in S-\eta$ the derived model at $\l$ satisfies MC (mouse capturing). Otherwise the derived model at such a $\lambda$ would have an initial segment satisfying $ZF+AD^++\theta_2=\Theta$ and the aforementioned result of Steel would still imply that there is a mouse $\M$ as above. Thus, we assume that\\\\

(1) for all $\l\in S-\eta$ the derived model at $\l$ satisfies $AD^++MC+\theta_1=\Theta$.\\\\

Fix $\l\in S-\eta$. Let $h\subseteq Coll(\omega, <\l)$ be $V[g]$-generic. Let $M$ be the derived model at $\l$ as computed by $h\cap Coll(\omega, <\l)$. Let $(\P, \Sigma)\in M$ be a hod pair such that 
\begin{center}
$M\models ``\l^\P=0$ and $\Sigma$ is a super fullness preserving iteration strategy with branch condensation that is strongly guided by $\vec{A}$" 
\end{center}
where $\vec{A}$ is some semi-scale on the universal $\Pi^2_1$ set of $M$ consisting of $OD^M$ prewellorderings. 

We have that in $V[g*h]$, $\Sigma\in (Hom^*)^{V[\mathbb{R}^*]}$. Let again $f(\l)<\l$ be the least cardinal $<\l$ such that it is forced by $Coll(\omega, <\l)$ that a $<\l$-uB code of $\Sigma$ appears after collapsing $f(\l)$. Again we obtain that $f$ is constant on a tail of $S$. Let $\nu$ be this constant value and let $k\subseteq Coll(\omega, \nu)$ be $V[g]$-generic. Then in $V[g][k]$, we have a pair $(\P, \Sigma)$ such that\\\\

(2) $\l^\P=0$, $\Sigma$ is universally Baire iteration strategy for $\P$ as witnessed by some class size trees $T, U$ with $p[T]=\Sigma$ such that for every $\l\in S-\nu$, if $h\subseteq Coll(\omega, <\l)$ is $V[g][k]$-generic, $M$ is the derived model at $\l$ computed by $h$ and $\Lambda=(p[T])^{V(\bR^*)}$ then $M\models ``\Lambda$ is a super fullness preserving iteration strategy with branch condensation that is strongly guided by $\vec{A}$, which is a semi-scale on the universal $\Pi^2_1$-set consisting of $OD$ prewellorderings". \\\\

Given a $V[g][k]$-generic $h$ we let $\Sigma^h=(p[T])^{V[g][k][h]}$. We now work in $W=V[g][k]$. We will abuse the notation and use $\Sigma$ instead of $\Sigma^{g*k}$. 

Suppose now that $D\subseteq Ord$ and $f:D\rightarrow Ord$ is an increasing function such that for $\a\in D$, $f[\a]\subseteq \a$, if $\a$ is not the maximum element of $D$ then $f(\a)\in D$ and $f(\a)$ is a strong limit cardinal in $W$. We then let $\Q_f$ be the output of the fully backgrounded $\mathcal{J}^{\vec{E}, \Sigma}$-construction done as follows. Let $(\xi_\a: \a\in \gg)$ be an enumeration of $D$ in increasing order. We also assume that $f(\xi_{\a+1})=\xi_{\a+2}$. Here we allow $\gg=Ord$. We then have the following conditions.
\begin{enumerate}
\item $\Q_f|f(\xi_0)$ is the output of $\mathcal{J}^{\vec{E}, \Sigma}$ construction done inside $W_{f(\xi_0)}$ using extenders with critical point $>\xi_0$.
\item $\Q_f|f(\xi_{\a+2})$ is the output of the $\mathcal{J}^{\vec{E}, \Sigma}$ construction done inside $W_{f(\xi_{\a+2})}$ over $\Q_f|f(\xi_\a)$ using extenders with critical point $>f(\xi_{\a+1})=\xi_{\a+2}$.
\item If $\l\in D$ is a limit ordinal then $\Q_f|\l=\cup_{\a<\l}\Q_f|\xi_\a$. 
\item If for any $\a<\b$, some level of $\Q_f|\xi_\b$ projects across $\xi_\a$ then we stop the construction. 
\end{enumerate}
Otherwise we let $\Q_f$ be the final model of the construction. 

 Notice that if there is a class size inner model with a Woodin cardinal that is a limit of Woodin cardinals then we are done proving our theorem. Otherwise, it follows from the results of \cite{Neeman}, that for every $f$, if $Ord\not \subseteq \Q_f$ then it is because for some $\a<\b$, some level of $\Q_f|\xi_\b$ projects across $\xi_\a$. 
 
Given a transitive set $X$, let $\W^\Sigma(X)$ be the union of all sound $\Sigma$-mice over $X$ that project to $X$ and are fully iterable. Also, let $Lp^{\Sigma, +}(X)$ be the union of all sound $\Sigma$-mice over $X$ that project to $X$ and whose countable submodels are fully iterable in all generic extensions. Recall that $Lp^\Sigma(X)$ is the union of all sound $\Sigma$-mice over $X$ that project to $X$ and whose countable submodels are $\omega_1+1$-iterable.\\

We now prove a claim that shows that the $Lp^\Sigma$ operator stabilizes at various derived model.\\

\textit{Claim 1.} There is a stationary class $S_0$ such that 
\begin{enumerate}
\item if $\l\in S_0$ then $\l$ is a limit of Woodin cardinals,
\item whenever $h\subseteq Coll(\omega, <\eta)$ is $W$-generic and $M$ is the derived model computed in $W(\bR^*)$ then letting $\Phi=\Sigma^h\rest HC^M$, for any $X\in HC^M$ and $\M\insegeq Lp^{\Phi}(X)$ such that $\rho_\omega(\M)=o(X)$, $W[h]\models ``\M$ has a universally Baire iteration strategy".
\end{enumerate}
\begin{proof} Suppose not. Then the set of $\l\in S$ for which clause 2 above fails is stationary. Given such a $\l$, we let $f(\l)<\l$ be least such that whenever $h\subseteq Coll(\omega, f(\l))$ is $W$-generic, there is a set of reals $A\in W[h]$ such that $A$ is $<\l$-universally Baire but not universally Baire. We then have that for a stationary class $E$, $f(\l)$ is constant. Let $\eta$ be this constant value and let $h\subseteq Coll(\omega, <\eta)$ be $W$-generic. Then for each $\l\in E$ we have a set of reals $A_\l\in W[h]$ such that in $W[h]$, $A_\l$ is $<\l$-uB but not uB. It follows that for a stationary class $F\subseteq E$ and for some set of reals $A\in W[h]$, $A=A_\l$ for all $\l\in F$.
\end{proof}

Notice that if $\l\in S_0$ then $Lp^{\Sigma, +}$ operator on sets in $W_\l$ is just the $Lp^\Sigma$-operator of the derived model. Our next claim connects $\emptyset$-iterability with background constructions.\\

\textit{Claim 2.} Suppose $\l\in S_0$ and $h\subseteq Coll(\omega, <\l)$ is $W$-generic. Let $M$ be the derived model computed in $W(\bR^*)$ and let $\Phi=\Sigma^h\rest HC^M$. Let $\Q$ be a $\Phi$-suitable premouse that is $\emptyset$-iterable. Let $\eta<\l$ be such that $\Q\in W_\eta[h\cap Coll(\omega, \eta)]$. Let $\nu\in (\eta, \l)$ be the least such that $Lp^{\Sigma, +}(W_\nu)\models ``\nu$ is a Woodin cardinal". Let $A\in \powerset(\nu)\cap \mathcal{J}(W_\nu)$ code $W_\nu$  and let $\R^-=(\mathcal{J}^{\vec{E}, \Sigma})^{W_\nu}$ in which extenders used have critical point $>\eta$ and cohere $A$. Let $\R=Lp_{\omega}^{\Sigma, +}(\R^-)$. Then $\R$ is a correct iterate of $\Q$, and hence, in $M$, $\R$ is $\emptyset$-iterable.\\\\
\begin{proof} It follows from clause 2 of Claim 1 and the universality of background constructions (see Lemma 2.11 of \cite{ATHM}) that $\Q$ indeed iterates to $\R^-$ in the sense there is a correctly guided tree $\T$ on $\Q$ such that $\R^-=\M(\T)$. 

It is then enough to show that $\R\models ``\nu$ is a Woodin cardinal". Suppose not. We have that $A$ and hence $W_\nu$ is generic for the extender algebra of $\R^-$ at $\nu$. It then follows that $Lp^{\Sigma, +}(W_\nu)$, using $\S$-constructions (see Section 3.8 of \cite{ATHM}), can be translated into $\Sigma$-mouse over $\R^-$. It therefore follows that $Lp^{\Sigma, +}(\R^-)\models ``\nu$ is a Woodin cardinal".
\end{proof}

Notice that a standard application of Fodor's lemma gives as a stationary class $S_1$ and ordinals $\eta^*<\nu^*$ such that fixing any set $A\in \powerset(\nu^*)$ that codes $W_{\nu^*}$ and letting $\R^-=(\mathcal{J}^{\vec{E}, \Sigma})^{W_{\nu^*}}$ in which extenders used have critical point $>\eta^*$ and cohere $A$, then for every $\l\in S_1$, $Lp_\omega^{\Sigma, +}(\R^-)$ is $\emptyset$-iterable in the derived model at $\l$. Let then $\P_0=Lp_\omega^{\Sigma, +}(\R^-)$.

Let now $k_0:S_1\rightarrow Ord$ be defined by letting $k_0(\a)$ be the least $W$-cardinal $>\a$ such that $Lp^{\Sigma, +}(W_{k_0(\a)})\models ``k_0(\a)$ is a Woodin cardinal" and for any $X\in W_\a$, letting $A\in \powerset(k_0(\a))\cap \mathcal{J}(W_{k_0(\a)})$ code $W_{k_0(\a)}$ and $\R^-$ be the output of $\mathcal{J}^{\vec{E}, \Sigma}$ construction of $W_{k_0(\a)}$ done over $X$ using extenders with critical point $>\a$, $Lp_\omega^{\Sigma, +}(\R^-)$ is an $\emptyset$-iterable $\Sigma$-suitable premouse in the derived model at any $\l\in S_1-\a$ (here we use the relativised version of Claim 2).

Suppose now that $D=(\xi_\a: \a<\gg)\subseteq S_1$ is such that for every $\a<\gg$, $k_0[\xi_\a]\subseteq \xi_\a$ and let
$f_D: D\cup rng(k_0\rest D)\rightarrow Ord$ be given by 
\begin{center}
$f_D(\xi)=\begin{cases}
k_0(\xi) &: \xi\in D\\
\eta &: \xi\not \in D \wedge \eta=min(D-\xi)\\
\end{cases}$
\end{center}

\textit{Claim 3.} Suppose $\l\in S_1$ and $\R\in W_\l$ are such that for some $\zeta$, $\R$ is a $\Sigma$-suitable premouse in the derived model at $\l$ and letting $\d$ be the largest Woodin cardinal of $\R$, $\d$ is a $W$-cardinal. Then there is $\eta\geq \l$ such that for any $D\subseteq Ord-\eta$ such that $f_D$ is defined, every level of $\Q_{f_D}$ construction build over $\R$ preserves the Woodiness of $\d$.  \\\\
\begin{proof}
Suppose not. We can then produce class many backgrounded constructions $(\Q_\a: \a\in Ord)$ and strictly increasing ordinals $(\eta_\a: \a\in Ord)$ such that for every $\a$, each $\Q_\a$ is done over $\R$ and size $\d$ hulls of every model appearing in the construction producing $\Q_\a$ is $\eta_\a$-iterable. It then follows that for some $\a$, for all $\b\geq \a$, $\Q_\b|(\d^+)^{\Q_\b}=\Q_\a|(\d^+)^{\Q_\a}$. Fix such $\a$, we get that $\Q_\a|(\d^+)^{\Q_\a}\insegeq Lp^{\Sigma, +}(\R|\d)$, contradiction.
\end{proof}

Given $\l\in S_1$, let $\eta_\l$ be the least ordinal that witnesses the validity of  the claim for any $\R\in W_\l$. We now define a set $D=(\l_\a: \a\in Ord)$ by induction as follows. 
\begin{enumerate}
\item $\l_0=min(S_1-o(\P_0))$,
\item letting $\l_{\a+1}^*=min(S_1-k_0(\l_\a))$, $\l_{\a+1}=min(S_1-\eta_{\l_{\a+1}^*})$,
\item for limit ordinals $\xi$, $\l_\xi=\sup_{\a<\xi}\l_\a$
\end{enumerate}

It now follows from Claim 3 that for any even $\a$, $k_0(\l_\a)$ is a Woodin cardinal in $\Q_{f_D}$ which is done over $\P_0$. This can be proven by a simple induction using Claim 3. \\\\

\textit{Claim 4.} $\Q_{f_D}$ is a translatable structure.\\\\
\begin{proof}
Notice that it follows from our construction that for any $\Q_{f_D}$-cutpoint $\xi>o(\P_0)$, if $\eta$ is the least Woodin cardinal of $\Q_{f_D}$ above $\xi$ and $\l=min(S_1-\xi)$ then $\Q_{f_D}|(\eta^{+\omega})^{\Q_{f_D}}$ is an $\emptyset$-iterable $\Sigma$-suitable premouse in the derived model at $\l$. This observation takes care of clause 3 of \rdef{translatable structure}. Clause 4 is satisfies as otherwise we already have the model that we are trying to build. The rest of the clause are consequence of our choice of $\Sigma$ and can be easily shown by using generic interpretability and generic comparisons (see Theorem 3.10 and Lemma 3.35 of \cite{ATHM})
\end{proof}

Applying \rthm{main theorem 1} we obtain a model as in the statement of \rthm{trev question}.
\end{proof}

Our final corollary shows that the two theories considered in \rthm{trev question} are equiconsistent.

\begin{corollary} The follows theories are equiconsistent.
\begin{enumerate}
\item $(T_1:)$ $ZFC+``$there is a stationary class of $\l$ such that $\l$ is a limit of Woodin cardinals and the derived model at $\l$ satisfies $\theta_0<\Theta$"
\item $(T_2:)$ $ZFC+``$there is a proper class of Woodin cardinals and a strong cardinal"
\end{enumerate}
\end{corollary}
\begin{proof}
It follows from \rthm{trev question} that $Con(T_1)\rightarrow Con(T_2)$. The reverse direction is a theorem of Woodin that can be found in \cite{DMT}. Woodin showed that under the hypothesis of clause 2, if $\l$ is a limit of Woodin cardinals bigger than a strong cardinal then the derived model at $\l$ satisfies $\theta_0<\Theta$.
\end{proof}

\bibliographystyle{plain}
\bibliography{TrevQues}

\begin{thebibliography}{10}

\bibitem{Neeman}
Itay Neeman.
\newblock Inner models in the region of a {W}oodin limit of {W}oodin cardinals.
\newblock {\em Ann. Pure Appl. Logic}, 116(1-3):67--155, 2002.

\bibitem{hodlsa}
Grigor Sargsyan.
\newblock The analysis of ${HOD}$ below the theory ${AD}^+$ +``the largest
  {S}muslin cardinal is a member of the {S}olovay sequence, available at
  math.rutgers.edu/$\sim$gs481.

\bibitem{complsa}
Grigor Sargsyan.
\newblock The comparison theory of hod pairs below ${AD}^+$ +``the largest
  {S}muslin cardinal is a member of the {S}olovay sequence, available at
  math.rutgers.edu/$\sim$gs481.

\bibitem{ATHM}
Grigor Sargsyan.
\newblock Hod mice and the {M}ouse {S}et {C}onjecture, available at
  http://math.rutgers.edu/$\sim$gs481/.

\bibitem{BSL}
Grigor Sargsyan.
\newblock Descriptive inner model theory.
\newblock {\em Bull. Symbolic Logic}, 19(1):1--55, 2013.

\bibitem{StrengthPFA1}
Grigor Sargsyan.
\newblock Nontame mouse from the failure of square at a singular strong limit
  cardinal.
\newblock {\em J. Math. Log.}, 14(1):1450003, 47, 2014.

\bibitem{TrangSarg}
Grigor Sargsyan and Nam Trang.
\newblock Non-tame mice from tame failures of the unique branch hypothesis.
\newblock {\em Canad. J. Math.}, 66(4):903--923, 2014.

\bibitem{DMATM}
John~R. Steel.
\newblock Derived models associated to mice.
\newblock In {\em Computational prospects of infinity. {P}art {I}.
  {T}utorials}, volume~14 of {\em Lect. Notes Ser. Inst. Math. Sci. Natl. Univ.
  Singap.}, pages 105--193. World Sci. Publ., Hackensack, NJ, 2008.

\bibitem{DMT}
John~R. Steel.
\newblock The derived model theorem.
\newblock In {\em Logic {C}olloquium 2006}, Lect. Notes Log., pages 280--327.
  Assoc. Symbol. Logic, Chicago, IL, 2009.

\bibitem{Trang}
Nam Trang.
\newblock H{OD} in natural models of {${AD}^+$}.
\newblock {\em Ann. Pure Appl. Logic}, 165(10):1533--1556, 2014.

\end{thebibliography}

\end{document}